# OBJECTIVE PRIORS FOR THE BIVARIATE NORMAL MODEL

By James O. Berger[1] and Dongchu Sun[2]

*Duke University and University of Missouri-Columbia*

Study of the bivariate normal distribution raises the full range of issues involving objective Bayesian inference, including the different types of objective priors (e.g., Jeffreys, invariant, reference, matching), the different modes of inference (e.g., Bayesian, frequentist, fiducial) and the criteria involved in deciding on optimal objective priors (e.g., ease of computation, frequentist performance, marginalization paradoxes). Summary recommendations as to optimal objective priors are made for a variety of inferences involving the bivariate normal distribution.

In the course of the investigation, a variety of surprising results were found, including the availability of objective priors that yield exact frequentist inferences for many functions of the bivariate normal parameters, including the correlation coefficient.

## 1. Introduction and prior distributions.

1.1. *Notation and problem statement.* The bivariate normal distribution of $(x_1, x_2)'$ has mean parameters $\boldsymbol{\mu} = (\mu_1, \mu_2)'$ and covariance matrix

$$\boldsymbol{\Sigma} = \begin{pmatrix} \sigma_1^2 & \rho\sigma_1\sigma_2 \\ \rho\sigma_1\sigma_2 & \sigma_2^2 \end{pmatrix},$$

where $\rho$ is the correlation between $x_1$ and $x_2$. The density is

$$\frac{1}{2\pi\sigma_1\sigma_2\sqrt{1-\rho^2}}$$
$$\times \exp\left\{ -\frac{\sigma_2^2(x_1-\mu_1)^2 + \sigma_1^2(x_2-\mu_2)^2 - 2\rho\sigma_1\sigma_2(x_1-\mu_1)(x_2-\mu_2)}{2\sigma_1^2\sigma_2^2(1-\rho^2)} \right\}.$$

[1]Supported in part by NSF Grant DMS-01-03265.

[2]Supported in part by NSF Grant SES-03-51523 and NIH Grant R01-MH071418.

The data consists of an independent random sample $\mathbf{X} = (\mathbf{x}_k = (x_{1k}, x_{2k}), k = 1, \dots, n)$ of size $n \geq 3$, for which the sufficient statistics are

$$(1) \quad \overline{\mathbf{x}} = \begin{pmatrix} \overline{x}_1 \\ \overline{x}_2 \end{pmatrix} \quad \text{and} \quad \mathbf{S} = \sum_{k=1}^{n} (\mathbf{x}_k - \overline{\mathbf{x}})(\mathbf{x}_k - \overline{\mathbf{x}})' = \begin{pmatrix} s_{11} & r\sqrt{s_{11}s_{22}} \\ r\sqrt{s_{11}s_{22}} & s_{22} \end{pmatrix},$$

where, for $i, j = 1, 2$,

$$\overline{x}_i = n^{-1} \sum_{j=1}^{n} x_{ij}, \qquad s_{ij} = \sum_{k=1}^{n} (x_{ik} - \overline{x}_i)(x_{jk} - \overline{x}_j) \quad \text{and} \quad r = \frac{s_{12}}{\sqrt{s_{11}s_{22}}}.$$

We will denote prior densities as $\pi(\mu_1, \mu_2, \sigma_1 \ \sigma_2, \ \rho)$, and the corresponding posterior densities as $\pi(\mu_1, \mu_2, \sigma_1 \ \sigma_2, \ \rho \mid \mathbf{X})$ (all with respect to $d\mu_1 \, d\mu_2 \, d\sigma_1 \, d\sigma_2 \, d\rho$).

We consider objective inference for parameters of the bivariate normal distribution and functions of these parameters, with special focus on development of objective confidence or credible sets. Section 1.2 introduces many of the key issues to be covered, through a summary of some of the most interesting results involving priors yielding exact frequentist procedures; this section also raises interesting historical and philosophical issues. For easy access, Section 1.3 presents our summary recommendations as to which priors to utilize.

Often, the posteriors for the recommended priors are essentially available in computational closed form, allowing direct Monte Carlo simulation. Section 2 provides simple accept-reject schemes for computing with the recommended priors in other cases. Sections 3 and 4 develop the needed theory, concerning what are called reference priors and matching priors, respectively, and also present various simulations that were conducted to enable summary recommendations to be made.

*Notation:* In addition to $(\mu_1, \mu_2, \sigma_1, \sigma_2, \rho)$, the following parameters will be considered:

$$(2) \quad \eta_1 = \frac{1}{\sigma_1}, \qquad \eta_2 = \frac{1}{\sigma_2\sqrt{1-\rho^2}}, \qquad \eta_3 = -\frac{\rho}{\sigma_1\sqrt{1-\rho^2}},$$

$$\theta_1 = \frac{\rho\sigma_2}{\sigma_1}, \qquad \theta_2 = \sigma_2^2(1-\rho^2), \qquad \theta_3 \equiv |\mathbf{\Sigma}| = \sigma_1^2\sigma_2^2(1-\rho^2),$$

$$(3)$$

$$\theta_4 = \frac{\sigma_2\sqrt{1-\rho^2}}{\sigma_1},$$

$$\theta_5 = \frac{\mu_1}{\sigma_1}, \qquad \theta_6 = \sigma_1^2\sigma_2^2, \qquad \theta_7 = \frac{\sigma_2}{\sigma_1}, \qquad \theta_8 = \frac{\mu_2}{\sigma_2},$$

$$(4)$$

$$\theta_9 \equiv \sigma_{12} = \rho\sigma_1\sigma_2,$$

$$(5) \quad \theta_{10} = \sigma_1^2 + \sigma_2^2 - 2\rho\sigma_1\sigma_2,$$

(6)     $\theta_{11} = \mathbf{d}'\boldsymbol{\Sigma}\mathbf{d}$     $[\mathbf{d}' = (d_1, d_2)$ not proportional to $(0,1)]$,

(7)     $\lambda_1 = ch_{\max}(\boldsymbol{\Sigma}), \qquad \lambda_2 = ch_{\min}(\boldsymbol{\Sigma}).$

Some of these parameters have straightforward statistical interpretations. Since $(x_2 \mid x_1, \boldsymbol{\mu}, \boldsymbol{\Sigma}) \sim N(\mu_2 + \theta_1(x_1 - \mu_1), \theta_2)$, it is clear that $\theta_1$ is a regression coefficient, $\theta_2$ is a conditional variance, and $\eta_2^2$ is the corresponding precision. For the marginal distribution of $x_1$, $\eta_1^2$ is the precision and $\theta_5$ is the reciprocal of the coefficient of variation. $\theta_3$ is usually called the generalized variance. $(\eta_1, \eta_2, \eta_3)$ gives a type of Cholesky decomposition of the precision matrix $\boldsymbol{\Sigma}^{-1}$ [see (13) in Section 2.1]. $\theta_{10}$ is the variance of $x_1 - x_2$, and $\theta_{11}$ is the variance of $d_1 x_1 + d_2 x_2$. Finally, $\lambda_1$ and $\lambda_2$ are the largest and smallest eigenvalues of $\boldsymbol{\Sigma}$.

*Technical issue.* We will assume that $|\rho| < 1$ and $|r| < 1$ in virtually all expressions and results that follow. This is because, if either equals 1 in absolute value, then $\rho = \{\text{sign of } r\}$ with probability 1 (either frequentist or Bayesian posterior, as relevant). Indeed, the situation then essentially collapses to the univariate version of the problem, which is standard.

1.2. *Matching, constructive posteriors and fiducial distributions.* The bivariate normal distribution has been extensively studied from frequentist, fiducial and objective Bayesian perspectives. Table 1 summarizes a number of interesting results.

- For a variety of parameters, it presents objective priors (discussed below) for which the resulting Bayesian posterior credible sets of level $1 - \alpha$ are also exact frequentist confidence sets at the same level; in this case, the priors are said to be *exact frequentist matching*. This is a very desirable situation: see [23] and [2] for general discussion and the many earlier references.
- For $\mu_1, \mu_2, \sigma_1, \sigma_2$ and $\rho$, the constructive posterior distributions are also the fiducial distributions for the parameters, as found in Fisher [14, 15] and [21].
- Posterior distributions are presented as *constructive random distributions*, that is, by a description of how to simulate from them. Thus to simulate from the posterior distribution of $\sigma_1$, given the data (actually, only $s_{11}$ is needed), one draws independent $\chi^2_{n-1}$ random variables and simply computes the corresponding $\sqrt{s_{11}/\chi^2_{n-1}}$; this yields an independent sample from the fiducial/posterior distribution of $\sigma_1$.

Table 1 also lists the objective prior distributions that yield the indicated objective posterior. The notation $\pi_{ab}$ in the table stands for the important

Table 1

*Parameters with exact matching priors of the form $\pi_{ab}$, and associated constructive posteriors: Here $Z^*$ is a standard normal random variable, and $\chi^2_{n-1}$ and $\chi^2_{n-2}$ are chi-squared random variables with the indicated degrees of freedom, all random variables being independent. For $\mu_1, \mu_2, \sigma_1, \sigma_2$ and $\rho$, the indicated posteriors are also fiducial distributions*

| Parameter | Prior | Posterior |
|---|---|---|
| $\mu_1$ | $\pi_{1b}, \forall b$ (including $\pi_J$ and $\pi_H$) | $\overline{x}_1 + \dfrac{Z^*}{\sqrt{\chi^2_{n-1}}}\sqrt{\dfrac{s_{11}}{n}}$ |
| $\mu_2$ | $\pi_J = \pi_{10}$ | $\overline{x}_2 + \dfrac{Z^*}{\sqrt{\chi^2_{n-1}}}\sqrt{\dfrac{s_{22}}{n}}$ |
| $\mathbf{d}'\binom{\mu_1}{\mu_2}, \mathbf{d} \in \mathbb{R}^2$ | $\pi_J = \pi_{10}$ and $\pi_{H*}$ (see Table 4) | $\mathbf{d}'(\overline{x}_1, \overline{x}_2)' + \dfrac{Z^*}{\sqrt{\chi^2_{n-1}}}\sqrt{\dfrac{\mathbf{d}'\mathbf{Sd}}{n}}$ |
| $\sigma_1$ | $\pi_{1b}, \forall b$ (including $\pi_J$ and $\pi_H$) | $\dfrac{\sqrt{s_{11}}}{\sqrt{\chi^2_{n-1}}}$ |
| $\rho$ | $\pi_H = \pi_{12}$ | $\psi\!\left(\dfrac{-Z^*}{\sqrt{\chi^2_{n-1}}} + \dfrac{\sqrt{\chi^2_{n-2}}}{\sqrt{\chi^2_{n-1}}}\dfrac{r}{\sqrt{1-r^2}}\right)$ $\psi(y) = y/\sqrt{1+y^2}$ |
| $\eta_3 = -\dfrac{\rho}{\sigma_1\sqrt{1-\rho^2}}$ | $\pi_{a2}, \forall a$ (including $\pi_H$) | $\dfrac{Z^*}{\sqrt{s_{11}}} - \dfrac{\sqrt{\chi^2_{n-2}}}{\sqrt{s_{11}}}\dfrac{r}{\sqrt{1-r^2}}$ |
| $\theta_1 = \dfrac{\rho\sigma_2}{\sigma_1}$ | $\pi_{a2}, \forall a$ (including $\pi_H$) | $\dfrac{r\sqrt{s_{22}}}{\sqrt{s_{11}}} - \dfrac{Z^*}{\sqrt{\chi^2_{n-2}}}\dfrac{\sqrt{1-r^2}\sqrt{s_{22}}}{\sqrt{s_{11}}}$ |
| $\theta_2 = \sigma_2^2(1-\rho^2)$ | $\pi_{a2}, \forall a$ (including $\pi_H$) | $\dfrac{s_{22}(1-r^2)}{\chi^2_{n-2}}$ |
| $\theta_3 = |\mathbf{\Sigma}|$ | $\pi_H = \pi_{12}$ and $\pi_{IJ} = \pi_{21}$ | $\dfrac{|\mathbf{S}|}{\chi^2_{n-1}\chi^2_{n-2}}$ |
| $\theta_4 = \dfrac{\sigma_2\sqrt{1-\rho^2}}{\sigma_1}$ | $\pi_H = \pi_{12}$ | $\dfrac{\sqrt{\chi^2_{n-1}}}{\sqrt{\chi^2_{n-2}}}\dfrac{\sqrt{s_{22}(1-r^2)}}{\sqrt{s_{11}}}$ |
| $\theta_5 = \dfrac{\mu_1}{\sigma_1}$ | $\pi_{1b}, \forall b$ (including $\pi_J$ and $\pi_H$) | $\dfrac{Z^*}{\sqrt{n}} + \dfrac{\overline{x}_1\sqrt{\chi^2_{n-1}}}{\sqrt{s_{11}}}$ |
| $\mathbf{d}'\mathbf{\Sigma d}$ | $\pi_J = \pi_{10}$ and $\pi_{H*}$ (see Table 4) | $\dfrac{\mathbf{d}'\mathbf{Sd}}{\chi^2_{n-1}}$ |

class of prior densities (a subclass of the *generalized Wishart distributions* of [8])

$$(8) \qquad \pi_{ab}(\mu_1, \mu_2, \sigma_1, \sigma_2, \rho) = \frac{1}{\sigma_1^{3-a}\sigma_2^{2-b}(1-\rho^2)^{2-b/2}}.$$

Special cases of this class are the *Jeffreys-rule* prior $\pi_J = \pi_{10}$, the *right-Haar* prior $\pi_H = \pi_{12}$, the *independence Jeffreys* prior $\pi_{IJ} = \pi_{21} = \sigma_1^{-1}\sigma_2^{-1}(1-\rho^2)^{-3/2}$ and $\pi_{RO}$ which has $a = b = 1$. The independence Jeffreys prior follows from using a constant prior for the means, and then the Jeffreys prior for the covariance matrix with means given.

We highlight the results about $\rho$ in Table 1 because they are interesting from practical, historical and philosophical perspectives. First, it does not seem to be known that the indicated prior for $\rho$ is exact frequentist matching (proved here in Theorem 2). Indeed, standard statistical software utilizes various approximations to arrive at frequentist confidence sets for $\rho$, missing the fact that a simple exact confidence set exists, even for $n = 3$. It was, of course, known that exact frequentist confidence procedures could be constructed (cf. Exercise 54, Chapter 6 of [18]), but explicit expressions do not seem to be available.

The historically interesting aspect of this posterior for $\rho$ is that it is also the fiducial distribution of $\rho$. Geisser and Cornfield [16] studied the question of whether the fiducial distribution of $\rho$ could be reproduced as an objective Bayesian posterior, and they concluded that this was most likely not possible. The strongest evidence for this arose from Brillinger [7], which used results from [19] and a difficult analytic argument to show that there does not exist a prior $\pi(\rho)$ such that the fiducial density of $\rho$ equals $f(r \mid \rho)\pi(\rho)$, where $f(r \mid \rho)$ is the density of $r$ given $\rho$. Since the fiducial distribution of $\rho$ only depends on $r$, it was certainly reasonable to speculate that if it were not possible to derive this distribution from the density of $r$ and a prior, then it would not be possible to do so in general. The above result, of course, shows that this speculation was incorrect.

The philosophically interesting aspect of this situation is that Brillinger's result does show that the fiducial/posterior distribution for $\rho$ provides another example of the marginalization paradox ([13]). This leads to an interesting philosophical conundrum of a type that we have not previously seen: a complete fiducial/objective Bayesian/frequentist unification can be obtained for inference about $\rho$, but only if violation of the marginalization paradox is accepted. We will shortly introduce a prior distribution that avoids the marginalization paradox for $\rho$, but which is not exactly frequentist matching. We know of no way to adjudicate between the competing goals of exact frequentist matching and avoidance of the marginalization paradox, and so will simply present both as possible objective Bayesian approaches. (Note that the same conundrum also arises for $\theta_5 = \mu_1/\sigma_1$; the exact frequentist matching prior results in a marginalization paradox, as shown in [24].) Some interesting examples of improper priors resulting in marginalization paradox can be found from Ghosh and Yang [17] and Datta and Ghosh [10, 11].

### 1.3. *Recommended priors.*

It is actually rare to have exact matching priors for parameters of interest. Also, one is often interested in very complex functions of parameters (e.g., predictive distributions) and/or joint distributions of parameters. For such problems it is important to have a general objective prior that seems to perform reasonably well for all quantities of interest. Furthermore, it is unappealing to many Bayesians to change the prior according to which parameter is declared to be of interest, and an objective prior that performs well overall is often sought.

TABLE 2
*Recommendations of objective priors for various parameters in the bivariate normal model:* □ *indicates that the posterior will not be exact frequentist matching. (For $\mu_2$ and parameters with $\sigma_1$ replaced by $\sigma_2$, use the right-Haar prior with the variances interchanged.)*

| Prior | Parameter |
|-------|-----------|
| $\pi_{R\rho}$ | $\boxed{\rho}$ , $\boxed{\frac{\sigma_1}{\sigma_2}}$ , $\boxed{\text{general use}}$ |
| $\pi_H$ | $\mu_1,\ \sigma_1,\ \rho,\ \eta_3,\ \frac{\rho\sigma_2}{\sigma_1},\ \sigma_2^2(1-\rho^2),\ |\boldsymbol{\Sigma}|,\ \frac{\sigma_2}{\sigma_1}\sqrt{1-\rho^2},\ \frac{\mu_1}{\sigma_1}$ |
| $\tilde{\pi}_H$ (see Table 4) | $\mathbf{d}'(\mu_1,\mu_2)',\ \mathbf{d}'\boldsymbol{\Sigma}\mathbf{d}$ |
| $\pi_{R\lambda}$ | $\boxed{ch_{\max}(\boldsymbol{\Sigma})}$ |
| $\pi_{R\sigma}$ | $\boxed{\sigma_{12}=\rho\sigma_1\sigma_2}$ |

The five priors we recommend for various purposes are $\pi_J$, $\pi_H$,

$$\text{(9)} \qquad \pi_{R\rho} \propto \frac{1}{\sigma_1\sigma_2(1-\rho^2)}, \qquad \pi_{R\sigma} \propto \frac{\sqrt{1+\rho^2}}{\sigma_1\sigma_2(1-\rho^2)}$$

and

$$\text{(10)} \qquad \pi_{R\lambda} \propto \frac{1}{\sigma_1\sigma_2(1-\rho^2)\sqrt{(\sigma_1/\sigma_2-\sigma_2/\sigma_1)^2+4\rho^2}}.$$

The first prior in (9) was developed in [20] and was studied extensively in [1], where it was shown to be a one-at-a-time reference prior (see Section 3). The second prior in (9) is new and is derived in Section 3. $\pi_{R\lambda}$ was developed as a one-at-a-time reference prior in [25].

With these definitions, we can make our summary recommendations. Table 2 gives the four objective priors that are recommended for use, and indicates for which parameters (or functions thereof) they are recommended. These recommendations are based on three criteria: (i) the degree of frequentist matching, discussed in Section 4; (ii) being a one-at-a-time reference prior, discussed in Section 3; and (iii) ease of computation. The rationale for each of the entries in the table, based on these criteria, is given in Section 4.5.

Another commonly used prior is the "scale prior," $\pi_S \propto (\sigma_1\sigma_2)^{-1}$. The motivation that is often given for this prior is that it is "standard" to use $\sigma_i^{-1}$ as the prior for a standard deviation $\sigma_i$, while $-1 < \rho < 1$ is on a bounded set and so one can use a constant prior in $\rho$. We do not recommend this prior, but do consider its performance in Section 4.5.

**2. Computation.** In this paper, a constant prior is always used for $(\mu_1, \mu_2)$, so that

$$\text{(11)} \qquad \left( \begin{pmatrix} \mu_1 \\ \mu_2 \end{pmatrix} \Big| \boldsymbol{\Sigma}, \mathbf{X} \right) \sim N_2\left( \begin{pmatrix} \overline{x}_1 \\ \overline{x}_2 \end{pmatrix}, n^{-1}\boldsymbol{\Sigma} \right).$$

Generation from this conditional posterior distribution is standard, so the challenge of simulation from the posterior distribution requires only sampling from $(\sigma_1, \sigma_2, \rho \mid \mathbf{X})$.

The marginal likelihood of $(\sigma_1, \sigma_2, \rho)$ satisfies

$$(12) \qquad L_1(\sigma_1, \sigma_2, \rho) \propto \frac{1}{|\mathbf{\Sigma}|^{(n-1)/2}} \exp\left(-\frac{1}{2} \operatorname{trace}(\mathbf{S}\mathbf{\Sigma}^{-1})\right).$$

It is immediate that, under the priors $\pi_J$ and $\pi_{IJ}$, the marginal posteriors of $\mathbf{\Sigma}$ are Inverse Wishart $(\mathbf{S}^{-1}, n)$ and Inverse Wishart $(\mathbf{S}^{-1}, n-1)$, respectively.

Berger, Strawderman and Tang [4] gave a Metropolis–Hastings algorithm to generate from $(\sigma_1, \sigma_2, \rho \mid \mathbf{X})$ based on the prior $\pi_{R\lambda}$. The following sections deal with the other priors we consider.

2.1. *Marginal posteriors of $(\sigma_1, \sigma_2, \rho)$ under $\pi_{R\rho}$, $\pi_{R\sigma}$, $\tilde{\pi}_{R\sigma}$, and $\pi_S$.* For these priors, an independent sample from $\pi(\sigma_1, \sigma_2, \rho \mid \mathbf{X})$ can be obtained by the following acceptance-rejection algorithm:

*Simulation step.* Generate $(\sigma_1, \sigma_2, \rho)$ from the independence Jeffreys posterior $\pi_{IJ}(\sigma_1, \sigma_2, \rho \mid \mathbf{X})$ [the Inverse Wishart $(\mathbf{S}^{-1}, n-1)$ distribution] and, independently, sample $u \sim \text{Uniform}(0, 1)$.
*Rejection step.* Suppose $M \equiv \sup_{(\sigma_1, \sigma_2, \rho)} \frac{\pi(\sigma_1, \sigma_2, \rho)}{\pi_{IJ}(\sigma_1, \sigma_2, \rho)} < \infty$. If $u \le \pi(\sigma_1, \sigma_2, \rho)/ [M\pi_{IJ}(\sigma_1, \sigma_2, \rho)]$, accept $(\sigma_1, \sigma_2, \rho)$; else, return to *Simulation step*.

For each of the priors listed in Table 3, the key ratio, $\pi/\pi_{IJ}$, is listed in the table, along with the upper bound $M$, the *Rejection step* and the resulting acceptance probability for $\rho = 0.80, 0.95, 0.99$. The rejection algorithm is quite efficient for sampling these posteriors. Indeed, for $\rho \approx 0$, the algorithms accept with probability near one and, even for large $|\rho|$, the acceptance probabilities are very reasonable for the priors $\pi_{R\rho}$, $\pi_{R\sigma}$, and $\tilde{\pi}_{R\sigma}$. For large $|\rho|$, the algorithm is less efficient for the posteriors under the prior $\pi_S$, but even these acceptance rates may well be fine in practice, given the simplicity of the algorithm.

2.2. *Computation under $\pi_{ab}$.* The most interesting prior of this form (besides the Jeffreys and independence Jeffreys priors) is the right-Haar prior $\pi_H$, although other priors such as $\pi_{11}$ arise as reference priors, and hence are potentially of interest. While Table 1 gave an explicit form for the most important marginal posteriors arising from priors of this form, it is of considerable interest that essentially closed form generation from the full posterior of any prior of this form is possible (see, e.g., [8]). This is briefly reviewed in this section, since the expressions for the resulting constructive posteriors are needed for later results on frequentist coverage.

TABLE 3
*Ratio $\pi/\pi_{IJ}$, upper bound $M$, rejection step and acceptance probability for*
*$\rho = 0.80, 0.95, 0.99$, when $\pi = \pi_{R\rho}$, $\pi_{R\sigma}$, $\tilde{\pi}_{R\sigma}$, $\pi_S$ and $\pi_{MS}$*

| | | **Bound** | | **Acceptance probability** | | |
|---|---|---|---|---|---|---|
| **Prior** | **Ratio $\frac{\pi}{\pi_{IJ}}$** | $M$ | **Rejection Step** | $\rho = 0.80$ | $\rho = 0.95$ | $\rho = 0.99$ |
| $\pi_{R\rho}$ | $\sqrt{1-\rho^2}$ | $1$ | $u \le \sqrt{1-\rho^2}$ | 0.6000 | 0.3122 | 0.1410 |
| $\pi_{R\sigma}$ | $\sqrt{1-\rho^4}$ | $1$ | $u \le \sqrt{1-\rho^4}$ | 0.7684 | 0.4307 | 0.1985 |
| $\tilde{\pi}_{R\sigma}$ | $\sqrt{\frac{1-\rho^2}{2-\rho^2}}$ | $\frac{1}{\sqrt{2}}$ | $u \le \sqrt{\frac{2(1-\rho^2)}{2-\rho^2}}$ | 0.7276 | 0.4215 | 0.1975 |
| $\pi_S$ | $(1-\rho^2)^{3/2}$ | $1$ | $u \le (1-\rho^2)^{3/2}$ | 0.2160 | 0.0304 | 0.0028 |

It is most convenient to work with the parameters $(\eta_1, \eta_2, \eta_3)$ given in (2). This parameterization gives a type of Cholesky decomposition of the precision matrix $\boldsymbol{\Sigma}^{-1}$,

$$(13) \qquad \boldsymbol{\Sigma}^{-1} = \begin{pmatrix} \eta_1 & \eta_3 \\ 0 & \eta_2 \end{pmatrix} \begin{pmatrix} \eta_1 & 0 \\ \eta_3 & \eta_2 \end{pmatrix},$$

which accounts for the simplicity of ensuing computations. Note that (2) is equivalent to

$$(14) \qquad \sigma_1 = \frac{1}{\eta_1}, \qquad \sigma_2 = \frac{\sqrt{\eta_1^2 + \eta_3^2}}{\eta_1 \eta_2}, \qquad \rho = -\frac{\eta_3}{\sqrt{\eta_1^2 + \eta_3^2}}.$$

The prior $\pi_{ab}$ of (8) for $(\mu_1, \mu_2, \sigma_1, \sigma_2, \rho)$ transforms to the extended conjugate class of priors for $(\mu_1, \mu_2, \eta_1, \eta_2, \eta_3)$, given by $\pi_{ab}(\mu_1, \mu_2, \eta_1, \eta_2, \eta_3) = \eta_1^{-a} \eta_2^{-b}$.

LEMMA 1. *Consider the prior $\pi_{ab}$.*

(a) *The marginal posterior of $\eta_3$ given $(\eta_1, \eta_2; \mathbf{X})$ is $N(-\eta_2 r \sqrt{s_{22}/s_{11}}, 1/s_{11})$.*

(b) *The marginal posterior distributions of $\eta_1$ and $\eta_2$ are independent and*

$$(\eta_1^2 \mid \mathbf{X}) \sim \text{Gamma}(\tfrac{1}{2}(n-a), \tfrac{1}{2}s_{11});$$

$$(\eta_2^2 \mid \mathbf{X}) \sim \text{Gamma}(\tfrac{1}{2}(n-b), \tfrac{1}{2}s_{22}(1-r^2)).$$

See [5] for a proof of this result. We next present the constructive posteriors of $(\eta_1, \eta_2, \eta_3)$, and from these derive the constructive posteriors of $(\mu_1, \mu_2, \sigma_1, \sigma_2, \rho)$ and other parameters. All results follow directly from Lemma 1 and (14).

In presenting the constructive posteriors, we will use a star to represent a random draw from the implied distribution; thus $\mu_1^*$ will represent a random draw from its posterior distribution, $Z_1^*, Z_2^*, Z_3^*$ will be independent draws from the standard normal distribution, and $\chi^{2*}_{n-a}$ and $\chi^{2*}_{n-b}$ will be independent draws from chi-squared distributions with the indicated degrees of freedom. Recall that these constructive posteriors are not only useful for simulation, but will be the key to proving exact frequentist matching results.

FACT 1. (a) *The constructive posterior of* $(\eta_1, \eta_2, \eta_3)$ *given* $\mathbf{X}$ *can be expressed as*

$$(15)\qquad \eta_1^* = \sqrt{\frac{\chi^{2*}_{n-a}}{s_{11}}}, \qquad \eta_2^* = \sqrt{\frac{\chi^{2*}_{n-b}}{s_{22}(1-r^2)}},$$

$$\eta_3^* = \frac{Z_3^*}{\sqrt{s_{11}}} - \frac{\sqrt{\chi^{2*}_{n-b}}}{\sqrt{s_{11}}}\frac{r}{\sqrt{1-r^2}}.$$

(b) *The constructive posterior of* $(\sigma_1, \sigma_2, \rho)$ *given* $\mathbf{X}$ *can be expressed as*

$$(16)\qquad \sigma_1^* = \sqrt{\frac{s_{11}}{\chi^{2*}_{n-a}}},$$

$$(17)\qquad \sigma_2^* = \sqrt{s_{22}(1-r^2)}\sqrt{\frac{1}{\chi^{2*}_{n-b}} + \frac{1}{\chi^{2*}_{n-a}}\left(\frac{Z_3^*}{\sqrt{\chi^{2*}_{n-b}}} - \frac{r}{\sqrt{1-r^2}}\right)^2},$$

$$(18)\qquad \rho^* = \psi(Y^*), \qquad Y^* = -\frac{Z_3^*}{\sqrt{\chi^{2*}_{n-a}}} + \frac{\sqrt{\chi^{2*}_{n-b}}}{\sqrt{\chi^{2*}_{n-a}}}\frac{r}{\sqrt{1-r^2}},$$

*where* $\psi(x) = x/\sqrt{1+x^2}$.

(c) *The constructive posterior for* $\mu_1$ *and* $\mu_2$ *can be written*

$$(19)\qquad \mu_1^* = \overline{x}_1 + \frac{Z_1^*}{\sqrt{\chi^{2*}_{n-a}}}\sqrt{\frac{s_{11}}{n}},$$

$$(20)\qquad \mu_2^* = \overline{x}_2 + \frac{Z_1^*}{\sqrt{\chi^{2*}_{n-a}}}\frac{r\sqrt{s_{22}}}{\sqrt{n}} + \left(\frac{Z_2^*}{\sqrt{\chi^{2*}_{n-b}}} - \frac{Z_3^*}{\sqrt{\chi^{2*}_{n-b}}}\frac{Z_1^*}{\sqrt{\chi^{2*}_{n-a}}}\right)\sqrt{\frac{s_{22}(1-r^2)}{n}}.$$

**3. Reference priors.** This paper began with an effort to derive and catalogue the possible reference priors for the bivariate normal distribution. The reference prior theory (cf. Bernardo [6] and Berger and Bernardo [3]) has arguably been the most successful technique for deriving objective priors. Reference priors depend on (i) specification of a parameter of interest; (ii) specification of nuisance parameters; (iii) specification of a grouping of parameters; and (iv) ordering of the groupings. These are all conveyed by the shorthand notation used in Table 4. Thus, $\{(\mu_1, \mu_2), (\sigma_1, \sigma_2, \rho)\}$ indicates that $(\mu_1, \mu_2)$ is the parameter of interest, with the others being nuisance parameters, and there are two groupings with the indicated ordering. (The resulting reference prior is the *independence Jeffreys* prior, $\pi_{IJ}$.) As another example, $\{\lambda_1, \lambda_2, \vartheta, \mu_1, \mu_2\}$ introduces the eigenvalues $\lambda_1 > \lambda_2$ of $\mathbf{\Sigma}$ as being primarily of interest, with $\vartheta$ (the angle defining the orthogonal matrix that diagonalizes $\mathbf{\Sigma}$), $\mu_1$ and $\mu_2$ being the nuisance parameters.

TABLE 4
*Reference priors for the bivariate normal model (where $\tilde{\mu}_1 = \mathbf{d}'(\mu_1, \mu_2)'$, $(\tilde{\sigma}_1)^2 = \theta_7$, $\tilde{\rho} = \mathbf{d}'\mathbf{\Sigma}(0,1)'/(\sigma_1\sqrt{\theta_7})$, $\tilde{\theta}_2 = \sigma_2^2[1-(\tilde{\rho})^2]$ and $\tilde{\theta}_1 = \tilde{\rho}\sigma_2/\tilde{\sigma}_1$); {{ }} indicates that any ordering of the parameters yields the same reference prior*

| Prior $\pi(\mu_1, \mu_2, \sigma_1, \sigma_2, \rho)$ | For parameter ordering | Has form (8) with |
|---|---|---|
| $\pi_J \propto \frac{1}{\sigma_1^2\sigma_2^2(1-\rho^2)^2}$ | $\{(\mu_1, \mu_2, \sigma_1, \sigma_2, \rho)\}$ | $(a,b) = (1,0)$ |
| $\pi_{IJ} \propto \frac{1}{\sigma_1\sigma_2(1-\rho^2)^{3/2}}$ | $\{(\mu_1, \mu_2), (\sigma_1, \sigma_2, \rho)\}$ | $(a,b) = (2,1)$ |
| $\pi_{R\rho} \propto \frac{1}{\sigma_1\sigma_2(1-\rho^2)}$ | $\{\rho, \sigma_1, \sigma_2\}, \{\theta_7, \theta_6, \rho\}$ | |
| $\pi_{R\sigma} \propto \frac{\sqrt{1+\rho^2}}{\sigma_1\sigma_2(1-\rho^2)}$ | $\{\sigma_1, \sigma_2, \rho\}$ | |
| $\tilde{\pi}_{R\sigma} \propto \frac{1}{\sigma_1\sigma_2(1-\rho^2)\sqrt{2-\rho^2}}$ | $\{\sigma_1, \rho, \sigma_2\}$ | |
| | $\{\sigma_1, \eta_3, \theta_2\}$ | |
| $\pi_{RO} \propto \frac{1}{\sigma_1^2\sigma_2(1-\rho^2)^{3/2}}$ | $\{\sigma_1, \theta_2, \eta_3\}$ | $(a,b) = (1,1)$ |
| $\pi_{R\lambda} \propto \frac{[((\sigma_1/\sigma_2)-(\sigma_2/\sigma_1))^2 + 4\rho^2]^{-1/2}}{\sigma_1\sigma_2(1-\rho^2)}$ | $\{\lambda_1, \lambda_2, \vartheta\}$ | |
| $\pi_H \propto \frac{1}{\sigma_1^2(1-\rho^2)}$ | $\{\{\sigma_1, \theta_1, \theta_2\}\}, \{\{\theta_1, \theta_3, \theta_4\}\}$ | $(a,b) = (1,2)$ |
| | $\{\{\eta_1, \eta_2, \theta_1\}\}, \{\{\eta_1, \theta_1, \theta_2\}\}$ | |
| $\tilde{\pi}_H \propto \frac{d\tilde{\mu}_1 d\mu_2 d\tilde{\sigma}_1 d\sigma_2 d\tilde{\rho}}{(\tilde{\sigma}_1)^2[1-(\tilde{\rho})^2]}$ | $\{\{\mathbf{d}'(\mu_1, \mu_2)', \mu_2, \theta_{11}, \tilde{\theta}_2, \tilde{\theta}_1\}\}$ | |

Based on experience with numerous examples, the reference priors that are typically judged to be best are one-at-a-time reference priors, in which each parameter is listed separately as its own group. Hence we will focus on these priors. It turns out to be the case that, for the one-at-a-time reference priors, the ordering of $\mu_1$ and $\mu_2$ among the variables is irrelevant. Hence if $\mu_1$ and $\mu_2$ are omitted from a listing in Table 4, the resulting reference prior is to be viewed as any one-at-a-time reference prior with the indicated ordering of other variables, with the $\mu_i$ being inserted anywhere in the ordering.

We are interested in finding one-at-a-time reference priors for the parameters $\mu_1, \mu_2, \sigma_1, \sigma_2, \rho, \eta_3, \theta_1, \ldots, \theta_9$ and $\lambda_1$. This is done in [5], with the results summarized in Table 4, for all these parameters (i.e., the parameter appears as the first entry in the parameter ordering) except $\eta_3$, $\sigma_{12}$, and $\mu_i/\sigma_i$; finding one-at-a-time reference priors for these parameters is technically challenging. (We do not explicitly list the reference priors for $\sigma_2$ in the table, since they can be found by simply switching with $\sigma_1$ in the various expressions.)

## 4. Comparisons of priors via frequentist matching.

4.1. *Frequentist coverage probabilities and exact matching.* Suppose a posterior distribution is used to create one-sided credible intervals $(\theta_L, \theta_{1-\alpha}(\mathbf{X}))$, where $\theta_L$ is the lower limit in the relevant parameter space and $\theta_{1-\alpha}(\mathbf{X})$ is the posterior quantile of the parameter $\theta$ of interest, defined by $P(\theta < \theta_{1-\alpha}(\mathbf{X}) \mid \mathbf{X}) = 1 - \alpha$. (Here $\theta$ is the random variable.) Of interest is the frequentist coverage of the corresponding confidence interval, that is, $C(\mu_1, \mu_2, \sigma_1, \sigma_2, \rho) = P(\theta < \theta_{1-\alpha}(\mathbf{X}) \mid \mu_1, \mu_2, \sigma_1, \sigma_2, \rho)$. (Here $\mathbf{X}$ is the random variable.) The closer $C(\mu_1, \mu_2, \sigma_1, \sigma_2, \rho)$ is to the nominal $1 - \alpha$, the better the procedure (and corresponding objective prior) is judged to be.

The main results about exact matching are given in Theorems 1 through 8. The proofs of Theorems 1, 2 and 8 are given in Section 5; the rest can be found in [5].

The following technical lemmas will be repeatedly utilized. The first lemma is from (3d.2.8) in [22]. Lemma 3 is easy.

LEMMA 2. *For $n \geq 3$ and given $\sigma_1, \sigma_2, \rho$, the following three random variables are independent and have the indicated distributions:*

$$(21) \quad T_2 = \left[ \frac{s_{11}}{\sigma_2^2(1-\rho^2)} \right]^{1/2} \left[ \frac{r\sqrt{s_{22}}}{\sqrt{s_{11}}} - \frac{\rho\sigma_2}{\sigma_1} \right] \equiv Z_3 \qquad (standard\ normal),$$

$$(22) \quad T_3 = \frac{s_{22}(1-r^2)}{\sigma_2^2(1-\rho^2)} \equiv \chi_{n-2}^2,$$

$$(23) \quad T_5 = \frac{s_{11}}{\sigma_1^2} \equiv \chi_{n-1}^2.$$

LEMMA 3. *Let $Y_{1-\alpha}$ denote the $1-\alpha$ quantile of any random variable $Y$.*

(a) *If $g(\cdot)$ is a monotonically increasing function, $[g(Y)]_{1-\alpha} = g(Y_{1-\alpha})$ for any $\alpha \in (0, 1)$.*

(b) *If $W$ is a positive random variable, $(WY)_{1-\alpha} \geq 0$ if and only if $Y_{1-\alpha} \geq 0$.*

We will reserve quantile notation for posterior quantiles, with respect to the $*$ distributions. Thus the quantile $[(\sigma_1 Z_3^* - r Z_3)/\chi_{n-1}^2 + \rho\sqrt{s_{11}}\chi_{n-b}^{2*}]_{1-\alpha}$ would be computed based on the joint distribution of $(Z_3^*, \chi_{n-b}^{2*})$, while holding $(\sigma_1, \rho, r, s_{11}, Z_3, \chi_{n-1}^2)$ fixed.

4.2. *Credible intervals for a class of functions of* $(\sigma_1, \sigma_2, \rho)$. We consider the one-sided credible intervals of $\sigma_1, \sigma_2$ and $\rho$ and some functions of the form

$$(24) \qquad \theta = \sigma_1^{d_1} \sigma_2^{d_2} g(\rho),$$

for $d_1, d_2 \in \mathbb{R}$ and some function $g(\cdot)$. We also consider a class of scale-invariant priors

$$(25) \qquad \pi(\mu_1, \mu_2, \sigma_1, \sigma_2, \rho) \propto \frac{h(\rho)}{\sigma_1^{c_1} \sigma_2^{c_2}},$$

for some $c_1, c_2 \in \mathbb{R}$ and a positive function $h$.

THEOREM 1. *Denote the* $1 - \alpha$ *posterior quantile of* $\theta$ *by* $\theta_{1-\alpha}(\mathbf{X})$ *under the prior (25). For any fixed* $(\mu_1, \mu_2, \sigma_1, \sigma_2, \rho)$, *the frequentist coverage of the credible interval* $(\theta_L, \theta_{1-\alpha}(\mathbf{X}))$ *depends only on* $\rho$. *Here* $\theta_L$ *is the lower boundary of the parameter space for* $\theta$.

Note that parameters $\rho, \eta_1, \eta_2, \eta_3, \theta_1, \ldots, \theta_4$ are all functions of the form (24). From Theorem 1, under any of the priors $\pi_J, \pi_{IJ}, \pi_{R\sigma}, \pi_{R\rho}, \pi_{RO}, \pi_H, \pi_S$, the frequentist coverage probabilities of credible intervals for any of these parameters will depend only on $\rho$. We will show that the frequentist coverage probabilities could be exact under the prior $\pi_{ab}$. Since $\eta_1(\eta_2)$ is a monotone function of $\sigma_1(\theta_2)$, we consider only $\rho$ and the last 5 parameters.

4.3. *Coverage probabilities under* $\pi_{ab}$.

THEOREM 2. (a) *For* $\psi$ *defined in (18), the posterior* $1 - \alpha$ *quantile of* $\rho$ *is* $\rho_{1-\alpha}^* = \psi(Y_{1-\alpha}^*)$. (b) *For any* $\alpha \in (0, 1)$, $\boldsymbol{\xi} = (\mu_1, \mu_2, \sigma_1, \sigma_2)$ *and* $\rho \in (-1, 1)$,

$$(26) \qquad \begin{aligned} &P(\rho < \rho_{1-\alpha}^* \mid \boldsymbol{\xi}, \rho) \\ &= P\left( \frac{\sqrt{1-\rho^2} Z_3 + \rho \sqrt{\chi_{n-1}^2}}{\sqrt{\chi_{n-2}^2}} > \left( \frac{\sqrt{1-\rho^2} Z_3^* + \rho \sqrt{\chi_{n-a}^{2*}}}{\sqrt{\chi_{n-b}^{2*}}} \right)_\alpha \mid \rho \right). \end{aligned}$$

(c) *(26) equals* $1 - \alpha$ *if and only if the right Haar prior is used, that is,* $(a, b) = (1, 2)$.

THEOREM 3. (a) *For any* $\alpha \in (0, 1)$, $\boldsymbol{\xi} = (\mu_1, \mu_2, \sigma_1, \sigma_2)$ *and* $\rho \in (-1, 1)$,

$$(27) \qquad \begin{aligned} &P(\eta_3 < (\eta_3^*)_{1-\alpha} \mid \boldsymbol{\xi}, \rho) \\ &= P\left( \frac{Z_3 + \frac{\rho}{\sqrt{1-\rho^2}} \sqrt{\chi_{n-1}^2}}{\sqrt{\chi_{n-2}^2}} < \left( \frac{Z_3^* + \frac{\rho}{\sqrt{1-\rho^2}} \sqrt{\chi_{n-1}^2}}{\sqrt{\chi_{n-b}^{2*}}} \right)_{1-\alpha} \mid \rho \right). \end{aligned}$$

(b) *(27) equals* $1 - \alpha$ *for any* $-1 < \rho < 1$ *if and only if* $b = 2$.

THEOREM 4. (a) *The constructive posterior of* $\theta_1 = \rho\sigma_2/\sigma_1$ *has the expression*

$$\theta_1^* = \frac{r\sqrt{s_{22}}}{\sqrt{s_{11}}} - \frac{Z_3^*}{\sqrt{\chi_{n-b}^{2*}}}\frac{\sqrt{1-r^2}\sqrt{s_{22}}}{\sqrt{s_{11}}}.$$

(b) *For any* $\alpha \in (0,1)$, $\boldsymbol{\xi} = (\mu_1, \mu_2, \sigma_1, \sigma_2)$ *and* $\rho \in (-1,1)$,

$$(28) \qquad P(\theta_1 < (\theta_1^*)_{1-\alpha} \mid \boldsymbol{\xi}, \rho) = P\left(t_{n-2} < \sqrt{\frac{n-2}{n-b}}(t_{n-b}^*)_{1-\alpha}\right),$$

*which does not depend on* $\rho$. *Furthermore,* *(28)* *equals* $1 - \alpha$ *if and only if* $b = 2$.

THEOREM 5. (a) *The constructive posterior of* $\theta_2 = \sigma_2^2(1-\rho^2)$ *is* $\theta_2^* = s_{22}(1-r^2)/\chi_{n-b}^{2*}$.

(b) *For any* $\alpha \in (0,1)$, $\boldsymbol{\xi} = (\mu_1, \mu_2, \sigma_1, \sigma_2)$ *and* $\rho \in (-1,1)$,

$$(29) \qquad P(\theta_2 < (\theta_2^*)_{1-\alpha} \mid \boldsymbol{\xi}, \rho) = P(\chi_{n-2}^2 > (\chi_{n-b}^{2*})_\alpha),$$

*which does not depend on* $\rho$. *Furthermore,* *(29)* *equals* $1 - \alpha$ *if and only if* $b = 2$.

THEOREM 6. (a) *The constructive posterior of* $\theta_3 = |\boldsymbol{\Sigma}|$ *is* $\theta_3^* = |\mathbf{S}|/(\chi_{n-a}^{2*}\chi_{n-b}^{2*})$.

(b) *For any* $\boldsymbol{\xi} = (\mu_1, \mu_2, \sigma_1, \sigma_2)$ *and* $\rho \in (-1,1)$,

$$(30) \qquad P(\theta_3 < (\theta_3^*)_{1-\alpha} \mid \boldsymbol{\xi}, \rho) = P(\chi_{n-1}^2 \chi_{n-2}^2 > (\chi_{n-a}^{2*}\chi_{n-b}^{2*})_\alpha),$$

*which does not depend on* $\rho$. *Furthermore,* *(30)* *equals* $1 - \alpha$ *iff* $(a,b)$ *is* $(1,2)$ *or* $(2,1)$.

THEOREM 7. (a) *The constructive posterior of* $\theta_4$ *is*

$$\theta_4^* = \sqrt{\frac{\chi_{n-a}^{2*}}{\chi_{n-b}^{2*}}}\sqrt{\frac{s_{22}(1-r^2)}{s_{11}}}.$$

(b) *For any* $\boldsymbol{\xi} = (\mu_1, \mu_2, \sigma_1, \sigma_2)$ *and* $\rho \in (-1,1)$,

$$(31) \qquad P(\theta_4 < (\theta_4^*)_{1-\alpha} \mid \boldsymbol{\xi}, \rho) = P(\chi_{n-1}^2/\chi_{n-2}^2 < (\chi_{n-a}^{2*}/\chi_{n-b}^{2*})_{1-\alpha}),$$

*which does not depend on* $\rho$. *Furthermore,* *(31)* *equals* $1 - \alpha$ *iff* $(a,b) = (1,2)$.

An interesting function of $(\mu_1, \mu_2, \sigma_1, \sigma_2, \rho)$ not of the form *(24)* is $\theta_5 = \mu_1/\sigma_1$.

THEOREM 8.  (a) *The constructive posterior of $\theta_5 = \mu_1/\sigma_1$ is*

$$\theta_5^* = \frac{Z_1^*}{\sqrt{n}} + \frac{\overline{x}_1}{\sqrt{s_{11}}}\sqrt{\chi_{n-a}^{2*}}.$$

(b) *For any $\alpha \in (0,1)$, the frequentist coverage of the credible interval $(-\infty, (\theta_5^*)_{1-\alpha})$ is*

$$(32) \qquad \begin{aligned} &P(\theta_5 < (\theta_5^*)_{1-\alpha} \mid \mu_1, \mu_2, \sigma_1, \sigma_2, \rho) \\ &= P\left(\frac{Z_1 - \theta_5\sqrt{n}}{\sqrt{\chi_{n-1}^2}} < \left(\frac{Z_1^* - \theta_5\sqrt{n}}{\sqrt{\chi_{n-a}^{2*}}}\right)_{1-\alpha} \,\Big|\, \theta_5\right), \end{aligned}$$

*which depends on $\theta_5$ only and equals $1 - \alpha$ if and only if $a = 1$.*

### 4.4. First order asymptotic matching.

Datta and Mukerjee [9] and Datta and Ghosh [12] discuss how to determine first-order matching priors for functions of parameters; these are priors such that the frequentist coverage of a one-sided credible interval is equal to the Bayesian coverage up to a term of order $n^{-1}$. For each of the nine objective priors $\pi_J$, $\pi_{IJ}$, $\pi_{R\rho}$, $\tilde{\pi}_{R\sigma}$, $\pi_{RO}$, $\pi_{R\lambda}$, $\pi_H$, $\pi_S$ and $\pi_{R\sigma}$, [5] determines if it is a first-order matching prior for each of the parameters $\mu_1, \mu_2, \sigma_1, \sigma_2, \rho, \eta_3, \theta_1, \ldots, \theta_{10}$. The results are listed in Table 5. For example, $\pi_J$ is a first order matching prior for $\mu_1, \mu_2, \sigma_1, \sigma_2$, $\theta_1$, $\theta_5$, $\theta_7$, $\theta_8$, and $\theta_{10}$, but not for $\eta_3$, $\theta_2$, $\theta_3$ and $\theta_9$.

### 4.5. Numerically computed coverage and recommendations.

First-order matching is only an asymptotic property, and finite sample performance is also crucial. We thus also implemented a modest numerical study, comparing the numerical values of frequentist coverages of the one-sided credible sets $P(\theta > q_{0.05})$ and $P(\theta < q_{0.95})$, for the parameters, $\theta$, listed in Table 6 and for the eight objective priors $\pi_J, \pi_{IJ}, \pi_{R\rho}, \pi_{R\sigma}, \pi_{RO}, \pi_{R\lambda}, \pi_H$ and $\pi_S$. As usual, $q_\alpha = q_\alpha(\mathbf{X})$ is the posterior $\alpha$-quantile of $\theta$, and the coverage probability is computed based on the sampling distribution of $q_\alpha(\mathbf{X})$ for the fixed parameter $(\mu_1, \mu_2, \sigma_1, \sigma_2)$ and $\rho$. Many of the coverage probabilities depend only on $\rho$, which was thus chosen to be the $x$-axis in the graphs. We considered the case $n = 3$ (the minimal possible sample size and hence the most challenging in terms of obtaining good coverage) and the two scenarios Case a: $(\mu_1, \mu_2, \sigma_1, \sigma_2) = (0, 0, 1, 1)$, and Case b: $(\mu_1, \mu_2, \sigma_1, \sigma_2) = (0, 0, 2, 1)$.

Here we present the numerical results concerning coverage for only two of the parameters: $\rho$ in Figure 1 and $\theta_7 = \sigma_2/\sigma_1$ in Figure 2. Table 6 summarizes the results from the entire numerical study, the details of which can be found in [5]. The recommendations made in Table 2 for the boxed parameters are justified from these numerical results as follows.

The inferences involving the nonboxed parameters in Table 2 are given in closed form in Table 1 (and so are computationally simple), and are exact

Table 5

*The first-order asymptotic matching of objective priors for $\mu_1, \mu_2, \sigma_1, \sigma_2, \rho$, $\mu_1 - \mu_2$, $\eta_3$, $\theta_j, j = 1, \ldots, 10$. Here a boldface letter indicates exact matching*

| | Asymptotic matching | |
|---|---|---|
| Prior $\pi(\boldsymbol{\mu_1}, \boldsymbol{\mu_2}, \boldsymbol{\sigma_1}, \boldsymbol{\sigma_2}, \rho)$ | **Yes** | **No** |
| $\pi_J \propto \frac{1}{\sigma_1^2 \sigma_2^2 (1-\rho^2)^2}$ | $\boldsymbol{\mu_1, \mu_2, \sigma_1, \sigma_2}$ | $\rho$ |
| | $\mu_1 - \mu_2, \theta_1, \boldsymbol{\theta_5}, \theta_7, \boldsymbol{\theta_8}, \theta_{10}$ | $\eta_3, \theta_2, \theta_3, \theta_9$ |
| $\pi_{IJ} \propto \frac{1}{\sigma_1 \sigma_2 (1-\rho^2)^{3/2}}$ | $\mu_1, \mu_2$ | $\sigma_1, \sigma_2, \rho$ |
| | $\mu_1 - \mu_2, \theta_1, \boldsymbol{\theta_3}, \theta_7$ | $\eta_3, \theta_2, \theta_5, \theta_8, \theta_9, \theta_{10}$ |
| $\pi_{R\rho} \propto \frac{1}{\sigma_1 \sigma_2 (1-\rho^2)}$ | $\mu_1, \mu_2, \rho$ | $\sigma_1, \sigma_2$ |
| | $\mu_1 - \mu_2, \theta_3, \theta_7$ | $\eta_3, \theta_1, \theta_2, \theta_5, \theta_8, \theta_9, \theta_{10}$ |
| $\tilde{\pi}_{R\sigma} \propto \frac{1}{\sigma_1 \sigma_2 (1-\rho^2)\sqrt{2-\rho^2}}$ | $\mu_1, \mu_2$ | $\sigma_1, \sigma_2, \rho$ |
| | $\mu_1 - \mu_2, \eta_3, \theta_7$ | $\theta_1, \theta_2, \theta_5, \theta_8, \theta_9, \theta_{10}$ |
| $\pi_{RO} \propto \frac{1}{\sigma_1^2 \sigma_2^2 (1-\rho^2)^{3/2}}$ | $\boldsymbol{\mu_1, \mu_2, \sigma_1}$ | $\sigma_2, \rho$ |
| | $\mu_1 - \mu_2, \theta_1, \boldsymbol{\theta_5}$ | $\eta_3, \theta_2, \theta_3, \theta_7, \theta_8, \theta_9, \theta_{10}$ |
| $\pi_{R\lambda} \propto \frac{[\sigma_1 \sigma_2 (1-\rho^2)]^{-1}}{\sqrt{((\sigma_1/\sigma_2)-(\sigma_2/\sigma_1))^2+4\rho^2}}$ | $\mu_1, \mu_2$ | $\sigma_1, \sigma_2, \rho$ |
| | $\mu_1 - \mu_2, \theta_3$ | $\eta_3, \theta_1, \theta_2, \theta_5,$ |
| | | $\theta_7, \theta_8, \theta_9, \theta_{10}$ |
| $\pi_H \propto \frac{1}{\sigma_1^4 (1-\rho^2)}$ | $\boldsymbol{\mu_1, \mu_2, \sigma_1, \rho}$ | $\sigma_2$ |
| | $\mu_1 - \mu_2, \boldsymbol{\eta_3, \theta_1, \theta_2, \theta_3, \theta_4, \theta_5}$ | $\theta_7, \theta_8, \theta_9, \theta_{10}$ |
| $\pi_S \propto \frac{1}{\sigma_1 \sigma_2}$ | $\mu_1, \mu_2$ | $\sigma_1, \sigma_2, \rho$ |
| | $\mu_1 - \mu_2, \theta_3, \theta_7$ | $\eta_3, \theta_1, \theta_2, \theta_5, \theta_8, \theta_9, \theta_{10}$ |
| $\pi_{R\sigma} \propto \frac{\sqrt{1+\rho^2}}{\sigma_1 \sigma_2 (1-\rho^2)}$ | $\mu_1, \mu_2$ | $\sigma_1, \sigma_2, \rho$ |
| | $\mu_1 - \mu_2, \theta_3, \theta_7, \theta_9$ | $\theta_1, \theta_2, \eta_3, \theta_5, \theta_8, \theta_{10}$ |

Table 6

*Performance of objective priors for each of the parameters*

| | Prior | | |
|---|---|---|---|
| Parameter | **Bad** | **Medium** | **Good** |
| $\mu_1$ | | rest | $\pi_{RO}, \pi_H, \pi_J$ |
| $\mu_1 - \mu_2$ | | rest | $\pi_J, \pi_{RO}$ |
| $\sigma_1$ | $\pi_{IJ}$ | rest | $\pi_H, \pi_{R\lambda}, \pi_{MS}$ |
| $\sigma_2$ | $\pi_H, \pi_{RO}, \pi_{IJ}$ | rest | $\pi_J$ |
| $\rho$ | $\pi_J, \pi_{IJ}, \pi_S, \pi_{RO}$ | | $\pi_{R\rho}, \pi_{R\sigma}, \pi_{R\lambda}, \pi_H, \pi_{MS}$ |
| $\lambda_1$ | rest | $\pi_J, \pi_{R\lambda}, \pi_{RO}$ | |
| $\theta_3 = \lvert \boldsymbol{\Sigma} \rvert$ | $\pi_{RO}, \pi_J$ | rest | $\pi_{IJ}, \pi_H$ |
| $\theta_7 = \frac{\sigma_2}{\sigma_1}$ | $\pi_H, \pi_J, \pi_{RO}, \pi_{R\lambda}$ | rest | |
| $\theta_9 = \sigma_{12}$ | $\pi_J, \pi_{IJ}$ (due to size) | rest | $\pi_H, \pi_{R\rho}, \pi_{R\sigma}$ |

frequentist matching. Furthermore, with the exception of $\mu_1/\sigma_1$ and $\eta_3$, the nonboxed parameters have the indicated priors as one-at-a-time reference priors, so all three criteria point to the indicated recommendation.

For $\rho$, we recommend using $\pi_{R\rho}$, since this prior is a one-at-a-time-reference for $\rho$, first-order matching (as shown in Table 5), and has excellent numerical coverage as shown in Figure 1. Note that some might prefer to use the right-Haar prior because of its exact matching for $\rho$ (even though it exhibits a marginalization paradox). For $\sigma_2/\sigma_1$, the one-at-a-time reference prior was also $\pi_{R\rho}$. As this was first-order frequentist matching and among the best in terms of numerical coverage (see Figure 2), we also recommend it for this parameter.

For $\lambda_1$, the situation is unclear. The one-at-a-time reference prior is $\pi_{R\lambda}$ and is hence our recommendation, but first-order matching results for this parameter are not known, and the numerical coverages of all priors were rather bad. For $\sigma_{12}$, the only first-order matching prior among our candidates is $\pi_{R\sigma}$. It also had the best numerical coverages, and so is a clear recommendation. Note, however, that we were not able to determine if it is a one-at-a-time reference prior for $\sigma_{12}$, so the recommendation should be considered tentative.

The most interesting question is what to recommend for general use, as an all-purpose prior. Looking at Table 2, it might seem that $\pi_H$ or even $\pi_J$ would be good choices, since they are optimal for so many parameters.

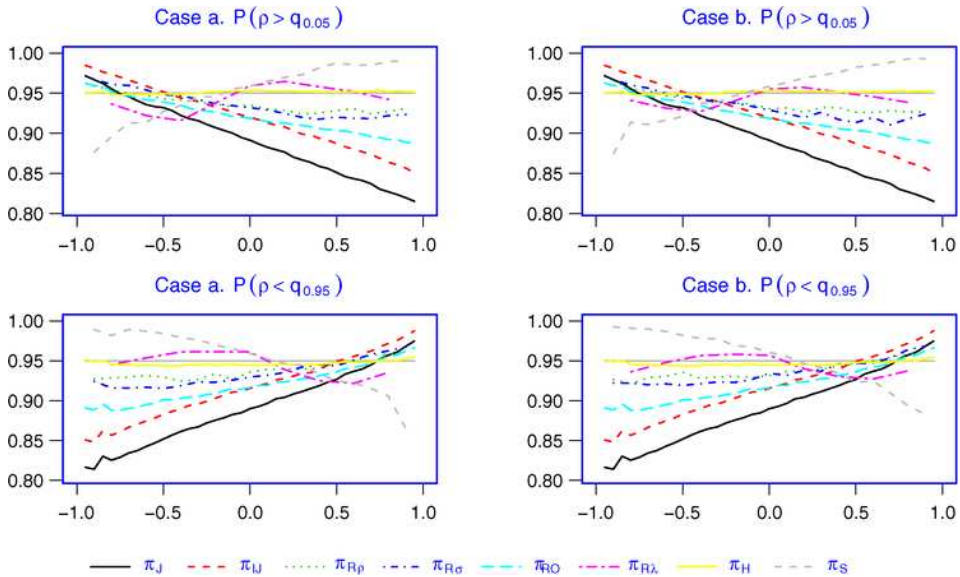

FIG. 1. *Frequentist coverages for $\rho$, where* Case a: $(\mu_1, \mu_2, \sigma_1, \sigma_2) = (0, 0, 1, 1)$, *and* Case b: $(\mu_1, \mu_2, \sigma_1, \sigma_2) = (0, 0, 2, 1)$. *The x-axis is for $\rho \in (-1, 1)$.*

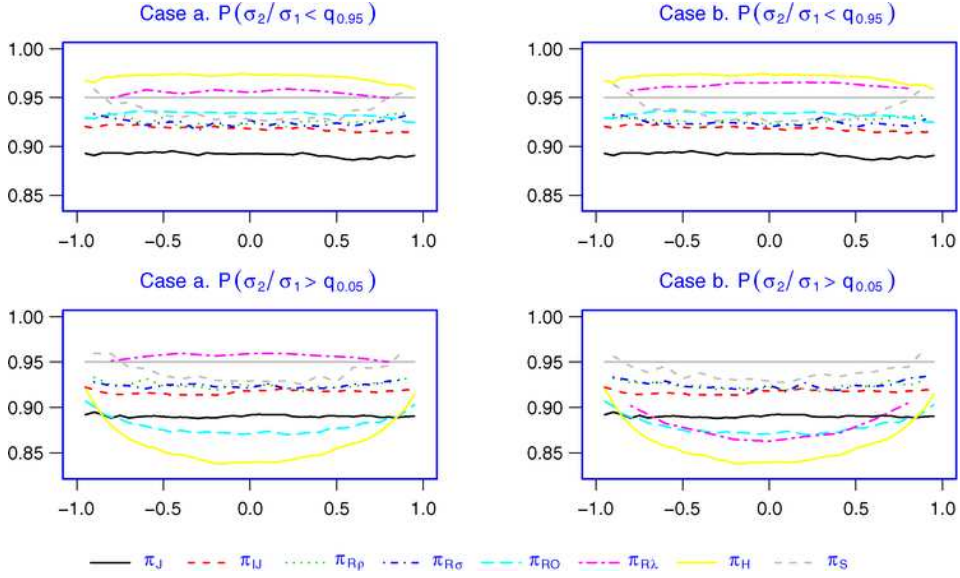

FIG. 2. *Frequentist coverages for* $\theta_7 = \sigma_2/\sigma_1$, *where* Case a: $(\mu_1, \mu_2, \sigma_1, \sigma_2) = (0, 0, 1, 1)$ *and* Case b: $(\mu_1, \mu_2, \sigma_1, \sigma_2) = (0, 0, 2, 1)$. *The x-axis is for* $\rho \in (-1, 1)$.

However, both these priors can also give quite bad coverages, as indicated in Figure 2 for $\pi_H$ and in Figures 1 and 2 for $\pi_J$. Indeed, from Table 6, the only priors that did not have significantly poor performance for at least one parameter (other than $\lambda_1$, for which no prior gave good coverages) were $\pi_{R\rho}$ and $\pi_{R\sigma}$. The numerical coverages for $\pi_{R\rho}$ and $\pi_{R\sigma}$ are virtually identical for all the parameters, so there is no principled way to choose between them. $\pi_{R\rho}$ is a commonly used prior and somewhat simpler, so it becomes our recommended choice for a general prior.

**5. Proofs.** Due to space limitations, we give only the proofs of Theorems 1, 2 and 8, because their proofs are quite different. The proofs of the other theorems in Section 4 are relatively easy consequences of Fact 1 and Lemmas 1–3. For details of these other proofs, see [5].

5.1. *Proof of Theorem* 1. With the constant prior for $(\mu_1, \mu_2)$, the marginal likelihood of $(\sigma_1, \sigma_2, \rho)$ depends on **S** and is proportional to

$$|\boldsymbol{\Sigma}|^{-(n-1)/2} \exp\{-\tfrac{1}{2}\operatorname{trace}(\mathbf{S}\boldsymbol{\Sigma}^{-1})\}.$$

Define

$$\mathcal{D} = \{(\sigma_1^*, \sigma_2^*, \rho^*) : \sigma_1^{*d_1} \sigma_2^{*d_2} g(\rho^*) < \sigma_1^{d_1} \sigma_2^{d_2} g(\rho)\},$$

$$G(\mathbf{X}, \sigma_1, \sigma_2, \rho) = \int_{\mathcal{D}} \pi(\sigma_1^*, \sigma_2^*, \rho^* \mid \mathbf{S}) \, d\sigma_1^* \, d\sigma_2^* \, d\rho^*.$$

Clearly, the frequentist coverage probability is

$$P\{\theta < \theta_{1-\alpha}(\mathbf{X}) \mid \mu_1, \mu_2, \sigma_1, \sigma_2, \rho\} = P\{G(\mathbf{S}, \sigma_1, \sigma_2, \rho) < 1 - \alpha \mid \sigma_1, \sigma_2, \rho\}.$$

Under the prior (25),

$$G(\mathbf{X}, \sigma_1, \sigma_2, \rho) = \frac{\int\int\int_{\mathcal{D}} \frac{h(\rho^*)\exp(-0.5\operatorname{trace}(\mathbf{S}\boldsymbol{\Sigma}^{*-1}))}{\sigma_1^{*(n-1+c_1)}\sigma_2^{*(n-1+c_2)}(1-\rho^{*2})^{(n-1)/2}}\,d\sigma_1^*\,d\sigma_2^*\,d\rho^*}{\int\int\int \frac{h(\rho^*)\exp(-0.5\operatorname{trace}(\mathbf{S}\boldsymbol{\Sigma}^{*-1}))}{\sigma_1^{*(n-1+c_1)}\sigma_2^{*(n-1+c_2)}(1-\rho^{*2})^{(n-1)/2}}\,d\sigma_1^*\,d\sigma_2^*\,d\rho^*},$$

where $\boldsymbol{\Sigma}^*$ is the $2 \times 2$ symmetric matrix, whose diagonal elements are $\sigma_1^{*2}$ and $\sigma_2^{*2}$, and off-diagonal element is $\sigma_1^*\sigma_2^*\rho^*$. Denote $\boldsymbol{\Xi} = \operatorname{diag}(1/\sigma_1, 1/\sigma_2)$ and make transformations

$$\mathbf{T} = \boldsymbol{\Xi}\mathbf{S}\boldsymbol{\Xi} = \begin{pmatrix} \dfrac{S_{11}}{\sigma_1^2} & \dfrac{S_{12}}{\sigma_1\sigma_2} \\[2mm] \dfrac{S_{12}}{\sigma_1\sigma_2} & \dfrac{S_{22}}{\sigma_2^2} \end{pmatrix} \quad \text{and} \quad \boldsymbol{\Omega} = \boldsymbol{\Xi}\boldsymbol{\Sigma}^*\boldsymbol{\Xi} = \begin{pmatrix} \omega_1^2 & \omega_1\omega_2\rho^* \\ \omega_1\omega_2\rho^* & \omega_2^2 \end{pmatrix}.$$

Clearly $\operatorname{trace}(\mathbf{S}\boldsymbol{\Sigma}^{*-1}) = \operatorname{trace}(\mathbf{T}\boldsymbol{\Omega}^{-1})$, and then

$$G(\mathbf{X}, \sigma_1, \sigma_2, \rho) = \frac{\int\int\int_{\widetilde{\mathcal{D}}} \frac{h(\rho^*)\exp(-0.5\operatorname{trace}(\mathbf{T}\boldsymbol{\Omega}^{-1}))}{\omega_1^{n-1+c_1}\omega_2^{n-1+c_2}(1-\rho^{*2})^{(n-1)/2}}\,d\omega_1\,d\omega_2\,d\rho^*}{\int\int\int \frac{h(\rho^*)\exp(-0.5\operatorname{trace}(\mathbf{T}\boldsymbol{\Omega}^{-1}))}{\omega_1^{n-1+c_1}\omega_2^{n-1+c_2}(1-\rho^{*2})^{(n-1)/2}}\,d\omega_1\,d\omega_2\,d\rho^*},$$

where $\widetilde{\mathcal{D}} = \{(\omega_1, \omega_2, \rho^*) : \omega_1^{d_1}\omega_2^{d_2}g(\rho^*) < g(\rho)\}$. Since the sampling distribution of $\mathbf{T}$ depends only on $\rho$, so does the sampling distribution of $G(\mathbf{X}, \sigma_1, \sigma_2, \rho)$. Also $\widetilde{\mathcal{D}}$ depends on $\rho$ only. The result thus holds.

5.2. *Proof of Theorem* 2. It follows from (18) and Lemma 3 (a) that

$$P(\rho < \rho_{1-\alpha}^* \mid \boldsymbol{\xi}, \rho) = P\left\{\left[\psi\left(\frac{-Z_3^*}{\sqrt{\chi_{n-a}^{2*}}} + \frac{\sqrt{\chi_{n-b}^{2*}}}{\sqrt{\chi_{n-a}^{2*}}}\frac{r}{\sqrt{1-r^2}}\right)\right]_{1-\alpha} > \rho \,\Big|\, \rho\right\},$$

Note that $\psi$, defined in (18), is invertible, and $\psi^{-1}(\rho) = \rho/\sqrt{1-\rho^2}$, for $|\rho| < 1$. It follows from Lemma 3 (a) and (b) that

$$P(\rho < \rho_{1-\alpha}^* \mid \boldsymbol{\xi}, \rho) = P\left(\left(\frac{-Z_3^*}{\sqrt{\chi_{n-a}^{2*}}} + \frac{\sqrt{\chi_{n-b}^{2*}}}{\sqrt{\chi_{n-a}^{2*}}}\frac{r}{\sqrt{1-r^2}} - \frac{\rho}{\sqrt{1-\rho^2}}\right)_{1-\alpha} > 0 \,\Big|\, \rho\right)$$

$$= P\left(\left(\frac{-Z_3^*}{\sqrt{\chi_{n-b}^{2*}}} - \frac{\rho}{\sqrt{1-\rho^2}}\frac{\sqrt{\chi_{n-a}^{2*}}}{\sqrt{\chi_{n-b}^{2*}}}\right)_{1-\alpha} + \frac{r}{\sqrt{1-r^2}} > 0 \,\Big|\, \rho\right).$$

It follows from (21)–(23) that

$$\frac{r}{\sqrt{1-r^2}} = \frac{s_{12}/\sqrt{s_{11}}}{\sqrt{s_{22}(1-r^2)}}$$

$$= \frac{\sigma_2\sqrt{1-\rho^2}Z_3 + (\rho\sigma_2/\sigma_1)\sqrt{s_{11}}}{\sigma_2\sqrt{1-\rho^2}\sqrt{\chi^2_{n-2}}}$$

$$= \frac{Z_3}{\sqrt{\chi^2_{n-2}}} + \frac{\rho}{\sqrt{1-\rho^2}}\frac{\sqrt{\chi^2_{n-1}}}{\sqrt{\chi^2_{n-2}}}.$$

Consequently,

$$P(\rho < \rho^*_{1-\alpha} \mid \boldsymbol{\xi}, \rho)$$

$$= P\Big(\frac{Z_3}{\sqrt{\chi^2_{n-2}}} + \frac{\rho}{\sqrt{1-\rho^2}}\frac{\sqrt{\chi^2_{n-1}}}{\sqrt{\chi^2_{n-2}}} < \Big(\frac{Z^*_3}{\sqrt{\chi^{2*}_{n-b}}} + \frac{\rho}{\sqrt{1-\rho^2}}\frac{\sqrt{\chi^{2*}_{n-a}}}{\sqrt{\chi^{2*}_{n-b}}}\Big)_{1-\alpha} \mid \rho\Big).$$

This completes the proof of part (a). For part (b), if (26) equals to $1-\alpha$ for any $-1 < \rho < 1$, choose $\rho = 0$ and get

$$P\Big(\frac{Z_3}{\sqrt{\chi^2_{n-2}}} < \Big(\frac{Z^*_3}{\sqrt{\chi^{2*}_{n-b}}}\Big)_{1-\alpha}\Big) = 1-\alpha,$$

which implies that $b = 2$. Substituting $b = 2$ into (26) shows that $a = 1$.

5.3. *Proof Theorem 8.* Part (a) is obvious. For part (b), since $\overline{x}_1 = \mu_1 + Z_1\sigma_1/\sqrt{n}$ and $Z_1$ and $\chi^2_{n-1}$ are independent, we have

$$(\theta_5 < (\theta^*_5)_{1-\alpha}) = \Big(\Big[\frac{Z^*_1}{\sqrt{n}} + \theta_5\Big(\sqrt{\frac{\chi^{2*}_{n-a}}{\chi^2_{n-1}}} - 1\Big) + \frac{Z_1}{\sqrt{n}}\sqrt{\frac{\chi^{2*}_{n-a}}{\chi^2_{n-1}}}\Big]_{1-\alpha} > 0\Big).$$

It follows from Lemma 3 (a) and (b) that

$$(\theta_5 < (\theta^*_5)_{1-\alpha}) = \Big(\Big[\frac{Z^*_1}{\sqrt{\chi^{2*}_{n-a}}} + \theta_5\Big(\frac{\sqrt{n}}{\sqrt{\chi^2_{n-1}}} - \frac{\sqrt{n}}{\sqrt{\chi^{2*}_{n-a}}}\Big) + \frac{Z_1}{\sqrt{\chi^2_{n-1}}}\Big]_{1-\alpha} > 0\Big)$$

$$= \Big(-\frac{Z_1}{\sqrt{\chi^2_{n-1}}} - \theta_5\frac{\sqrt{n}}{\sqrt{\chi^2_{n-1}}} < \Big(\frac{Z^*_1}{\sqrt{\chi^{2*}_{n-a}}} - \theta_5\frac{\sqrt{n}}{\sqrt{\chi^{2*}_{n-a}}}\Big)_{1-\alpha}\Big).$$

Because $Z_1$ and $-Z_1$ have the same distribution and $Z_1$ and $\chi^2_{n-1}$ are independent, (32) holds. If (32) equals $1-\alpha$ for any $\theta_5$, choose $\theta_5 = 0$,

$$P\Big(\frac{Z_1}{\sqrt{\chi^2_{n-1}}} < \Big(\frac{Z^*_1}{\sqrt{\chi^{2*}_{n-a}}}\Big)_{1-\alpha}\Big) = 1-\alpha,$$

which implies that $a = 1$. The result holds.

**Acknowledgments.** The authors are grateful to Fei Liu for performing the numerical frequentist coverage computations, to Xiaoyan Lin for computing the matching priors in Table 5, and to Susie Bayarri for helpful discussions. The authors gratefully acknowledge the very constructive comments of the editor, an associate editor and two referees.

ISDS
DUKE UNIVERSITY
BOX 90251
DURHAM, NORTH CAROLINA NC 27708-0251
USA
E-MAIL: berger@stat.duke.edu
URL: www.stat.duke.edu/˜berger

DEPARTMENT OF STATISTICS
UNIVERSITY OF MISSOURI-COLUMBIA
146 MIDDLEBUSH HALL
COLUMBIA, MISSOURI 65211-6100
USA
E-MAIL: sund@missouri.edu
URL: www.stat.missouri.edu/˜dsun